\input amstex
\input amsppt.sty
\magnification=\magstep1
\hsize=30truecc
\vsize=22.2truecm
\baselineskip=16truept
\NoBlackBoxes
\TagsOnRight \pageno=1 \nologo
\def\Z{\Bbb Z}
\def\N{\Bbb N}

\def\l{\left}
\def\r{\right}
\def\bg{\bigg}
\def\({\bg(}
\def\[{\bg\lfloor}
\def\){\bg)}
\def\]{\bg\rfloor}
\def\t{\text}
\def\f{\frac}

\def\bi{\binom}
\def\eq{\equiv}

\def\ls{\leqslant}
\def\gs{\geqslant}
\def\mo{\roman{mod}}

\def\al{\alpha}
\def\da{\delta}

\def\Proof{\noindent{\it Proof}}

\def\Remark{\medskip\noindent{\it  Remark}}

\hbox {Preprint, {\tt arXiv:0912.1280}}
\bigskip
\topmatter
\title Congruences involving binomial coefficients and Lucas sequences\endtitle
\author Zhi-Wei Sun\endauthor
\leftheadtext{Zhi-Wei Sun} \rightheadtext{Binomial coefficients and Lucas sequences}
\affil Department of Mathematics, Nanjing University\\
 Nanjing 210093, People's Republic of China
  \\  zwsun\@nju.edu.cn
  \\ {\tt http://math.nju.edu.cn/$\sim$zwsun}
\endaffil
\abstract In this paper we obtain some congruences involving central binomial coefficients and Lucas sequences.
For example, we show that if $p>5$ is a prime  then
$$\sum_{k=0}^{p-1}\f{F_k}{12^k}\bi{2k}k\eq\cases0\ (\mo\ p)&\t{if}\  p\eq\pm1\ (\mo\ 5),
\\1\ (\mo\ p)&\t{if}\ p\eq \pm13\ (\mo\ 30),\\-1\ (\mo\ p)&\t{if}\ p\eq\pm7\ (\mo\ 30),\endcases$$
where $\{F_n\}_{n\gs0}$ is the Fibonacci sequence.
We also raise several conjectures.
\endabstract
\thanks 2010 {\it Mathematics Subject Classification}.\,Primary 11B65;
Secondary 05A10,\,11A07, 11B39, 11E25.
\newline\indent {\it Keywords}. Central binomial coefficients, Lucas sequences, congruences modulo prime powers.
\newline\indent Supported by the National Natural Science
Foundation (grant 10871087) and the Overseas Cooperation Fund (grant 10928101) of China.
\endthanks
\endtopmatter
\document

\heading{1. Introduction}\endheading

Let $p$ be an odd prime. In 2003 Roderiguez-Villeags [RV] conjectured  that
$$\sum_{k=0}^{p-1}\f{\bi{2k}k^2}{16^k}\eq(-1)^{(p-1)/2}\ (\mo\ p^2)$$
and
$$\sum_{k=0}^{p-1}\f{\bi{2k}k^3}{64^k}\eq a(p)\ (\mo\ p^2),$$
where the sequence $\{a(n)\}_{n\gs1}$ is defined by
$$\sum_{n=1}^\infty a(n)q^n=q\prod_{n=1}^\infty(1-q^{4n})^6.$$
This was later confirmed by E. Mortenson [M1, M2] via
 the $p$-adic $\Gamma$-function and the Gross-Koblitz formula.
The reader may also consult [M3] and Ono [O] for more such ``super" congruences.

In a series of recent papers, the author [S09a-S09e] investigated congruences related to central binomial congruences
by using recurrences and combinatorial identities. (See also [PS] and [ST1, ST2].)

Let $A,B\in\Z$. The Lucas sequences $u_n=u_n(A,B)\ (n\in\N)$ and  $v_n=v_n(A,B)\ (n\in\N)$
are defined by
$$u_0=0,\ u_1=1,\ \t{and}\ u_{n+1}=Au_n-Bu_{n-1}\ (n=1,2,3,\ldots)$$
and
$$v_0=2,\ v_1=A,\ \t{and}\ v_{n+1}=Av_n-Bv_{n-1}\ (n=1,2,3,\ldots).$$
The characteristic equation $x^2-Ax+B=0$ has two roots
$$\al=\f{A+\sqrt{\Delta}}2\quad\t{and}\quad\beta=\f{A+\sqrt{\Delta}}2,$$
where $\Delta=A^2-4B$. It is well known that for any $n\in\N$ we have
$$u_n=\sum_{0\ls k<n}\al^k\beta^{n-1-k}\quad\t{and }\quad v_n=\al^n+\beta^n.$$
Note that $F_n=u_n(1,-1)$ and $L_n=v_n(1,-1)$ are Fibonacci numbers and Lucas numbers respectively.
The sequences $P_n=u_n(2,-1)$ and $Q_n=v_n(2,-1)$ are called the Pell sequence and its companion.
We also set $S_n=u_n(4,1)$ and $T_n=v_n(4,1)$ for $n\in\N$; the sequences $\{S_n\}_{n\gs0}$
and its companion $\{T_n\}_{n\gs0}$ are also useful (see, e.g., [S02]).

In this paper we study congruences involving both central binomial coefficients and Lucas sequences.
Now we state our main results.

\proclaim{Theorem 1.1}
Let $A,m\in\Z$ and let $p$ be an odd prime not dividing $m$.
Suppose that $\delta^2\eq A^2-4m^2\not\eq0\ (\mo\ p)$ where $\da\in\Z$.
Let $a,h\in\Z^+$. If $(\f{A+\da}{p^a})=(\f{2m}{p^a})$, then
$$\sum_{k=0}^{p^a-1}\f{u_k(A,m^2)\bi{2k}k^h}{m^k(-4)^{hk}}\eq0\ (\mo\ p).$$
If $(\f{A+\da}{p^a})=-(\f{2m}{p^a})$, then
$$\sum_{k=0}^{p^a-1}\f{v_k(A,m^2)\bi{2k}k^h}{m^k(-4)^{hk}}\eq0\ (\mo\ p).$$
\endproclaim
\proclaim{Corollary 1.1} Let $p\eq1\ (\mo\ 3)$ be a prime and let $a\in\Z^+$. Then
$$\sum_{k=0}^{p^a-1}\l(\f k3\r)\f{\bi{2k}k}{(-4)^k}
\eq\sum_{k=0}^{p^a-1}\l(\f k3\r)\f{\bi{2k}k^2}{16^k}
\eq\sum_{k=0}^{p^a-1}\l(\f k3\r)\f{\bi{2k}k^3}{(-64)^k}\eq0\ (\mo\ p).$$
When $p^a\eq 1\ (\mo\ 12)$, we have
$$\sum_{k=0}^{p^a-1}\l(\f k3\r)\f{\bi{2k}k}{4^k}\eq\sum_{k=0}^{p^a-1}\l(\f k3\r)\f{\bi{2k}k^2}{(-16)^k}
\eq\sum_{k=0}^{p^a-1}\l(\f k3\r)\f{\bi{2k}k^3}{64^k}\eq0\ (\mo\ p).$$
If $p^a\eq 7\ (\mo\ 12)$, then
$$\sum_{k=0}^{(p^a-1)/3}\f{\bi{6k}{3k}}{64^k}\eq0\ (\mo\ p)\ \t{and}\
\sum_{k=0}^{(p^a-1)/3}(-1)^k\f{\bi{6k}{3k}^2}{2^{12k}}\eq0\ (\mo\ p).$$
\endproclaim
\proclaim{Corollary 1.2} Let $p\eq\pm1\ (\mo\ 5)$ be a prime and let $a\in\Z^+$.
Then
$$\sum_{k=0}^{p^a-1}F_{2k}\f{\bi{2k}k}{(-4)^k}
\eq \sum_{k=0}^{p^a-1}F_{2k}\f{\bi{2k}k^2}{16^k}
\eq \sum_{k=0}^{p^a-1}F_{2k}\f{\bi{2k}k^3}{(-64)^k}\eq0\ (\mo\ p).$$
If $p^a\eq 1,9\ (\mo\ 20)$, then
$$\sum_{k=0}^{p^a-1}F_{2k}\f{\bi{2k}k}{4^k}\eq\sum_{k=0}^{p^a-1}F_{2k}\f{\bi{2k}k^2}{(-16)^k}
\eq \sum_{k=0}^{p^a-1}F_{2k}\f{\bi{2k}k^3}{64^k}\eq0\ (\mo\ p).$$
If $p^a\eq 11,19\ (\mo\ 20)$, then
$$\sum_{k=0}^{p^a-1}L_{2k}\f{\bi{2k}k}{4^k}\eq\sum_{k=0}^{p^a-1}L_{2k}\f{\bi{2k}k^2}{(-16)^k}
\eq \sum_{k=0}^{p^a-1}L_{2k}\f{\bi{2k}k^3}{64^k}\eq0\ (\mo\ p).$$
\endproclaim

\proclaim{Theorem 1.2} Let $p$ be an odd prime and let $A,B\in\Z$
and $p\nmid AB\Delta$, where $\Delta=A^2-4B$.

{\rm (i)} If $p\eq1\ (\mo\ 4)$, then
$$\sum_{k=0}^{p-1}\f{u_k(A,B)}{(16A)^k}\bi{2k}k^2\eq0\ (\mo\ p^2).$$
If $p\eq3\ (\mo\ 4)$, then
$$\sum_{k=0}^{p-1}\f{v_k(A,B)}{(16A)^k}\bi{2k}k^2\eq0\ (\mo\ p^2).$$

{\rm (ii)} Suppose that $(\f{\Delta}p)=1$. If
$(\f{-B}p)=1$, then
$$\sum_{k=0}^{p-1}\f{A^ku_k(A,B)}{(16B)^k}\bi{2k}k^2\eq0\ (\mo\ p).$$
If $(\f{-B}p)=-1$, then
$$\sum_{k=0}^{p-1}\f{A^kv_k(A,B)}{(16B)^k}\bi{2k}k^2\eq0\ (\mo\ p).$$
\endproclaim
\Remark. Theorem 1.2(i) with $A=1$ and $B=-1$ was first noted by R. Tauraso [T].

\proclaim{Corollary 1.3} Let $p\eq1\ (\mo\ 4)$ be a prime. Then
$$\sum_{k=0}^{p-1}\l(\f k3\r)\f{\bi{2k}k^2}{(-16)^k}\eq0\ (\mo\ p^2).$$
\endproclaim

\proclaim{Corollary 1.4} Let $p>5$ be a prime.

{\rm (i)} If $p\eq1,4\ (\mo\ 5)$ then
$$\sum_{k=0}^{p-1}\f{F_k}{(-16)^k}\bi{2k}k^2\eq0\ (\mo\ p).$$

{\rm (ii)} If $p\eq 1\ (\mo\ 4)$, then
$$\sum_{k=0}^{p-1}\f{F_{2k}}{48^k}\bi{2k}k^2\eq0\ (\mo\ p^2).$$
If $p\eq3\ (\mo\ 4)$ then
$$\sum_{k=0}^{p-1}\f{L_{2k}}{48^k}\bi{2k}k^2\eq0\ (\mo\ p^2).$$

{\rm (iii)} If $p\eq1,9\ (\mo\ 20)$, then
$$\sum_{k=0}^{p-1}\f{3^kF_{2k}}{16^k}\bi{2k}k^2\eq0\ (\mo\ p).$$
If $p\eq11,19\ (\mo\ 20)$, then
$$\sum_{k=0}^{p-1}\f{3^kL_{2k}}{16^k}\bi{2k}k^2\eq0\ (\mo\ p).$$
\endproclaim

\proclaim{Corollary 1.5} Let $p$ be an odd prime. If $p\eq 1\ (\mo\ 4)$, then
$$\sum_{k=0}^{p-1}\f{P_{k}}{32^k}\bi{2k}k^2\eq0\ (\mo\ p^2).$$
If $p\eq3\ (\mo\ 4)$ then
$$\sum_{k=0}^{p-1}\f{Q_{k}}{32^k}\bi{2k}k^2\eq0\ (\mo\ p^2).$$
If $p\eq\pm1\ (\mo\ 8)$, then
$$\sum_{k=0}^{p-1}\f{P_k}{(-8)^k}\bi{2k}k^2\eq0\ (\mo\ p).$$
\endproclaim

\proclaim{Corollary 1.6} Let $p>3$ be a prime. If $p\eq 1\ (\mo\ 4)$, then
$$\sum_{k=0}^{p-1}\f{S_{k}}{64^k}\bi{2k}k^2\eq0\ (\mo\ p^2).$$
If $p\eq3\ (\mo\ 4)$ then
$$\sum_{k=0}^{p-1}\f{T_{k}}{64^k}\bi{2k}k^2\eq0\ (\mo\ p^2).$$
If $p\eq1\ (\mo\ 12)$, then
$$\sum_{k=0}^{p-1}\f{S_k}{4^k}\bi{2k}k^2\eq0\ (\mo\ p).$$
If $p\eq11\ (\mo\ 12)$, then
$$\sum_{k=0}^{p-1}\f{T_k}{4^k}\bi{2k}k^2\eq0\ (\mo\ p).$$
\endproclaim

\proclaim{Theorem 1.3} Let $A,B\in\Z$ and $\Delta=A^2-4B$. Let $p$ be an odd prime and let $m\in\Z$ with $p\nmid m$.
Suppose that $p\nmid \Delta$ and $d^2\eq m^2-4Am+16B\not\eq0\ (\mo\ p)$ where $d\in\Z$.
Then
$$\sum_{k=0}^{p-1}\f{u_k(A,B)}{m^k}\bi{2k}k\eq\cases0\ (\mo\ p)&\t{if}\ (\f{\Delta}p)=1,
\\ -\f 4d(\f{2m}p)(\f{m-d-2A}p)\ (\mo\ p)
&\t{if}\ (\f{\Delta}p)=-1.\endcases$$
Also,
$$\sum_{k=0}^{p-1}\f{v_k(A,B)}{m^k}\bi{2k}k\eq\cases2(\f{2m}p)(\f{m-d-2A}p)\ (\mo\ p)&\t{if}\ (\f{\Delta}p)=1,
\\ \f {4A-2m}d(\f{2m}p)(\f{m-d-2A}p)\ (\mo\ p)
&\t{if}\ (\f{\Delta}p)=-1.\endcases$$
\endproclaim

Theorem 1.3 in the case $A=-1,\ B=1, \ m=-4$ and $d=1$ gives the following consequence.

\proclaim{Corollary 1.7} Let $p$ be an odd prime. Then
$$\sum_{k=0}^{p-1}\f{(\f k3)}{(-4)^k}\bi{2k}k\eq\f{(\f{-1}p)-(\f 3p)}2\ (\mo\ p).$$
\endproclaim

Applying Theorem 1.3 with $A=1,\ B=-1,\ m\in\{4,8\}$ and $d=4$, we immediately get the following corollary.

\proclaim{Corollary 1.8} Let $p$ be an odd prime. Then
$$\align\sum_{k=0}^{p-1}\f{F_k}{(-4)^k}\bi{2k}k\eq&\f{1-(\f p5)}2\ (\mo\ p),
\\ \sum_{k=0}^{p-1}\f{L_k}{(-4)^k}\bi{2k}k\eq&\f{5(\f p5)-1}2\ (\mo\ p),
\\\sum_{k=0}^{p-1}\f{F_k}{8^k}\bi{2k}k\eq&\l(\f 2p\r)\f{(\f p5)-1}2\ (\mo\ p),
\\\sum_{k=0}^{p-1}\f{L_k}{8^k}\bi{2k}k\eq&\l(\f 2p\r)\f{5(\f p5)-1}2\ (\mo\ p).
\endalign$$
\endproclaim

Theorem 1.3 in the case $A=2,\ B=-1,\ m\in\{-2,10\}$ and $d=2$, yields the following result.
\proclaim{Corollary 1.9} Let $p$ be an odd prime. Then
$$\align\sum_{k=0}^{p-1}\f{P_k}{(-2)^k}\bi{2k}k\eq&1-\l(\f 2p\r)\ (\mo\ p),
\\ \sum_{k=0}^{p-1}\f{Q_k}{(-2)^k}\bi{2k}k\eq&4\l(\f 2p\r)-2\ (\mo\ p).
\endalign$$
If $p\not=5$, then
$$\sum_{k=0}^{p-1}\f{P_k}{10^k}\bi{2k}k\eq\l(\f p5\r)\l(\l(\f2p\r)-1\r)\ (\mo\ p)$$
and
$$\sum_{k=0}^{p-1}\f{Q_k}{10^k}\bi{2k}k\eq\l(\f p5\r)\l(4\l(\f 2p\r)-2\r)\ (\mo\ p).$$
\endproclaim

Theorem 1.3 in the case $A=4,\ B=1,\ m\in\{1,15,16\}$, leads the following corollary.
\proclaim{Corollary 1.10} Let $p>3$ be a prime. Then
$$\align\sum_{k=0}^{p-1}S_k\bi{2k}k\eq&2\l(\l(\f p3\r)-\l(\f{-1}p\r)\r)\ (\mo\ p),
\\\sum_{k=0}^{p-1}T_k\bi{2k}k\eq&8\l(\f{-1}p\r)-6\l(\f p3\r)\ (\mo\ p),
\\\sum_{k=0}^{p-1}\f{S_k}{16^k}\bi{2k}k\eq&\f{(\f6p)-(\f 2p)}2\ (\mo\ p),
\\\sum_{k=0}^{p-1}\f{T_k}{16^k}\bi{2k}k\eq&3\l(\f6p\r)-\l(\f 2p\r)\ (\mo\ p).
\endalign$$
When $p>5$, we also have
$$\sum_{k=0}^{p-1}\f{S_k}{15^k}\bi{2k}k\eq 2\l(\f p5\r)\l(\l(\f 3p\r)-1\r)\ (\mo\ p)$$
and
$$\sum_{k=0}^{p-1}\f{T_k}{15^k}\bi{2k}k\eq 2\l(\f p5\r)\l(8\l(\f 3p\r)-6\r)\ (\mo\ p).$$
\endproclaim

\proclaim{Theorem 1.4} Let $A,B\in\Z$ and $\Delta=A^2-4B$. Let $p$ be an odd prime with $p\nmid AB$ and $(\f{\Delta}p)=1$.
Then
$$\sum_{k=0}^{p-1}\f{A^kv_k(A,B)}{(4B)^k}\bi{2k}k\eq2\l(\f{-B}p\r)\ (\mo\ p^2).$$
\endproclaim

\proclaim{Theorem 1.5} Let $p>5$ be a prime. Then
$$\sum_{k=0}^{p-1}\f{F_k}{12^k}\bi{2k}k\eq\cases0\ (\mo\ p)&\t{if}\  p\eq\pm1\ (\mo\ 5),
\\1\ (\mo\ p)&\t{if}\ p\eq \pm13\ (\mo\ 30),\\-1\ (\mo\ p)&\t{if}\ p\eq\pm7\ (\mo\ 30).\endcases$$
Also,
$$\sum_{k=0}^{p-1}\f{L_k}{12^k}\bi{2k}k\eq\cases-1\ (\mo\ p)&\t{if}\  p\eq\pm7\ (\mo\ 30),
\\1\ (\mo\ p)&\t{if}\ p\eq \pm13\ (\mo\ 30),\\2\ (\mo\ p)&\t{if}\ p\eq\pm1\ (\mo\ 30),
\\-2\ (\mo\ p)&\t{if}\ p\eq\pm11\ (\mo\ 30).\endcases$$
\endproclaim

Let $p>5$ be a prime. By the method we prove Theorem 1.5, we can also determine the following sums modulo $p$.
$$\sum_{k=0}^{p-1}\f{F_k}{m^k}\bi{2k}k\ \t{and}\ \sum_{k=0}^{p-1}\f{L_k}{m^k}\bi{2k}k\ (m=-3,\ 7,\ -8),$$
$$\sum_{k=0}^{p-1}\f{P_k}{m^k}\bi{2k}k\ \t{and}\ \sum_{k=0}^{p-1}\f{Q_k}{m^k}\bi{2k}k\ (m=-4,\ 12)$$
For example,
$$\sum_{k=0}^{p-1}\f{F_k}{(-3)^k}\bi{2k}k\eq\l(\f p5\r)-1\ (\mo\ p),$$
and
$$\sum_{k=0}^{p-1}\f{L_k}{(-3)^k}\bi{2k}k\eq\cases2\ (\mo\ p)&\t{if}\ p\eq\pm1\ (\mo\ 30),
\\-2\ (\mo\ p)&\t{if}\ p\eq\pm11\ (\mo\ 30),\\4\ (\mo\ p)&\t{if}\ p\eq\pm7\ (\mo\ 30),
\\-4\ (\mo\ p)&\t{if}\ p\eq\pm13\ (\mo\ 30).\endcases$$
Modifying the method slightly, we can also prove the following congruences.
$$\sum_{k=0}^{p-1}\f{C_kP_k}{(-2)^k}\eq2\l(\f 2p\r)-2\ (\mo\ p),\ \sum_{k=0}^{(p-1)/2}C_kS_k\eq\f{(\f p3)-1}2\ (\mo\ p),$$
$$\sum_{k=0}^{(p-1)/2}\f{C_kF_k}{(-4)^k}\eq 2\l(\f p5\r)-2\ (\mo\ p),$$
and
$$\sum_{k=0}^{(p-1)/2}\f{C_kF_k}{12^k}\eq\cases0\ (\mo\ p)&\t{if}\ p\eq\pm1\ (\mo\ 30),
\\4\ (\mo\ p)&\t{if}\ p\eq\pm7\ (\mo\ 30),\\8\ (\mo\ p)&\t{if}\ p\eq\pm13\ (\mo\ 30),
\\12\ (\mo\ p)&\t{if}\ p\eq\pm11\ (\mo\ 30).\endcases$$

\medskip

In the next section we will prove Theorem 1.1 and Corollaries 1.1 and 1.2.
Theorems 1.2-1.3 and Corollaries 1.3-1.6 will be proved in Section 3.
Section 4 is devoted to the proof of Theorems 1.4 and 1.5.
We are going to raise some challenging conjectures in Section 5.

\heading{2. Proofs of Theorem 1.1 and Corollaries 1.1-1.2}\endheading

\proclaim{Lemma 2.1}  Let $A,B\in\Z$ and let $\Delta=A^2-4B$. Suppose that $p$ is an odd prime and
$\da^2\eq\Delta\ (\mo\ p)$ with $\da\in\Z$. Then, for any $n\in\N$ we have
$$\da u_n(A,B)\eq\l(\f{A+\da}2\r)^n-\l(\f{A-\delta}2\r)^n\ (\mo\ p)\tag 2.1$$
and
$$v_n(A,B)\eq\l(\f{A+\da}2\r)^n+\l(\f{A-\delta}2\r)^n\ (\mo\ p).\tag2.2$$
\endproclaim
\Proof. Observe that
$$\align v_n(A,B)=&\l(\f{A+\sqrt{\Delta}}2\r)^n+\l(\f{A-\sqrt{\Delta}}2\r)^n
\\=&\f1{2^n}\sum_{k=0}^n\bi nkA^{n-k}\l(\sqrt{\Delta}^k+(-\sqrt{\Delta})^k\r)
=\f2{2^n}\sum^n\Sb k=0\\2\mid k\endSb\bi nkA^{n-k}\Delta^{k/2}
\\\eq&\f2{2^n}\sum^n\Sb k=0\\2\mid k\endSb\bi nkA^{n-k}\delta^k=\l(\f{A+\da}2\r)^n+\l(\f{A-\delta}2\r)^n\ (\mo\ p).
\endalign$$
If $p\mid\Delta$, then both sides of (2.1) are multiples of $p$. When $p\nmid\Delta$,
we have
$$\align u_n(A,B)=&\f1{\sqrt{\Delta}}\(\l(\f{A+\sqrt{\Delta}}2\r)^n-\l(\f{A-\sqrt{\Delta}}2\r)^n\)
\\=&\f1{2^n\sqrt{\Delta}}\sum_{k=0}^n\bi nkA^{n-k}\l(\sqrt{\Delta}^k-(-\sqrt{\Delta})^k\r)
\\=&\f2{2^n}\sum^n\Sb k=0\\2\nmid k\endSb\bi nkA^{n-k}\Delta^{(k-1)/2}
\\\eq&\f2{2^n}\sum^n\Sb k=0\\2\nmid k\endSb\bi nkA^{n-k}\delta^{k-1}
\\\eq&\f1{\delta}\(\l(\f{A+\da}2\r)^n-\l(\f{A-\delta}2\r)^n\)\ (\mo\ p).
\endalign$$
Thus both (2.1) and (2.2) hold. \qed

\proclaim{Lemma 2.2} Let $p$ be an odd prime and let $a\in\Z^+$.
Then,for every $k=0,\ldots,p^a-1$ we have
$$\bi{(p^a-1)/2}k\eq\f{\bi{2k}k}{(-4)^k}\ (\mo\ p).$$
\endproclaim
\Proof. The congruence appeared as [S09e, (2.3)]. \qed

\proclaim{Theorem 2.1} Let $A,B\in\Z$ and let $\Delta=A^2-4B$.
Let $p$ be an odd prime with $(\f{\Delta}p)=1$.
Suppose that $\delta^2\eq\Delta\not\eq0\ (\mo\ p)$ where $\da\in\Z$.
Let $a,h\in\Z^+$ and $m\in\Z$ with $m\not\eq0\ (\mo\ p)$. If $(\f B{p^a})=1$, then
$$\sum_{k=0}^{p^a-1}\f{u_k(A,B)}{m^k}\cdot\f{\bi{2k}k^h}{(-4)^{hk}}
\eq-\l(\f{2m(A+\da)}{p^a}\r)\sum_{k=0}^{p^a-1}\f{m^ku_k(A,B)}{B^k}\cdot\f{\bi{2k}k^h}{(-4)^{hk}}\ (\mo\ p)$$
and
$$\sum_{k=0}^{p^a-1}\f{v_k(A,B)}{m^k}\cdot\f{\bi{2k}k^h}{(-4)^{hk}}
\eq\l(\f{2m(A+\da)}{p^a}\r)\sum_{k=0}^{p^a-1}\f{m^kv_k(A,B)}{B^k}\cdot\f{\bi{2k}k^h}{(-4)^{hk}}\ (\mo\ p).$$
If $(\f B{p^a})=-1$, then
$$\sum_{k=0}^{p^a-1}\f{u_k(A,B)}{m^k}\cdot\f{\bi{2k}k^h}{(-4)^{hk}}
\eq\f1{\delta}\l(\f{2m(A+\da)}{p^a}\r)\sum_{k=0}^{p^a-1}\f{m^kv_k(A,B)}{B^k}\cdot\f{\bi{2k}k^h}{(-4)^{hk}}\ (\mo\ p)$$
and
$$\sum_{k=0}^{p^a-1}\f{v_k(A,B)}{m^k}\cdot\f{\bi{2k}k^h}{(-4)^{hk}}
\eq-\delta\l(\f{2m(A+\da)}{p^a}\r)\sum_{k=0}^{p^a-1}\f{m^ku_k(A,B)}{B^k}\cdot\f{\bi{2k}k^h}{(-4)^{hk}}\ (\mo\ p).$$
\endproclaim
\Proof. Set
$$n=\f{p^a-1}2,\ \ \al=\f{A+\da}2\ \ \t{and}\ \ \beta=\f{A-\da}2.$$
Clearly,
$$(2\al)^{(p-1)/2}=(A+\da)^{(p-1)/2}\eq\l(\f{A+\delta}p\r)\ (\mo\ p)$$
and hence
$$\al^n=\l(\al^{(p-1)/2}\r)^{\sum_{r=0}^{a-1}p^r}\eq\l(\f{2(A+\da)}p\r)^{\sum_{r=0}^{a-1}p^r}=\l(\f{2(A+\da)}{p^a}\r)\ (\mo\ p).$$
Similarly,
$$\beta^n\eq\l(\f{2(A-\da)}{p^a}\r)\ (\mo\ p)$$
and hence
$$\al^n\beta^n\eq\l(\f{A^2-\da^2}{p^a}\r)=\l(\f{4B}{p^a}\r)=\l(\f B{p^a}\r)\ (\mo\ p).$$

By Lemmas 2.1 and 2.2,
$$\align &\sum_{k=0}^{p-1}\f{u_k(A,B)}{m^k}\cdot\f{\bi{2k}k^h}{(-4)^{hk}}
\\\eq&\sum_{k=0}^n\bi nk^h\f{\al^k-\beta^k}{m^k\delta}=\sum_{k=0}^n\bi nk^h\f{\al^{n-k}-\beta^{n-k}}{m^{n-k}\delta}
\\\eq&\f1{m^n\da}\sum_{k=0}^n\bi nk^hm^k(\al^{n-k}-\beta^{n-k})
\\\eq&\f{(\f m{p^a})}{\da}\sum_{k=0}^n\bi nk^h\f{m^k}{B^k}(\al^n\beta^k-\beta^n\al^k)
\\\eq&\l(\f m{p^a}\r)\l(\f{2(A+\da)}p\r)\sum_{k=0}^n\f{\bi{2k}k^h}{(-4)^{hk}}\cdot\f{m^k}{B^k}\cdot\f{\beta^k-(\f B{p^a})\al^k}{\da}\ (\mo\ p).
\endalign$$
Similarly,
$$\align &\sum_{k=0}^{p-1}\f{v_k(A,B)}{m^k}\cdot\f{\bi{2k}k^h}{(-4)^{hk}}
\\\eq&\sum_{k=0}^n\bi nk^h\f{\al^k+\beta^k}{m^k}=\sum_{k=0}^n\bi nk^h\f{\al^{n-k}+\beta^{n-k}}{m^{n-k}}
\\\eq&\f1{m^n}\sum_{k=0}^n\bi nk^hm^k(\al^{n-k}+\beta^{n-k})
\\\eq&\l(\f m{p^a}\r)\sum_{k=0}^n\bi nk^h\f{m^k}{B^k}(\al^n\beta^k+\beta^n\al^k)
\\\eq&\l(\f m{p^a}\r)\l(\f{2(A+\da)}{p^a}\r)\sum_{k=0}^n\f{\bi{2k}k^h}{(-4)^{hk}}\cdot\f{m^k}{B^k}\l(\beta^k
+\l(\f B{p^a}\r)\al^k\r)\ (\mo\ p).
\endalign$$
Note that
$$\f{\beta^k-(\f B{p^a})\al^k}{\da}\eq\cases (\beta^k-\al^k)/\da\eq-u_k(A,B)\ (\mo\ p)&\t{if}\ (\f B{p^a})=1,
\\(\al^k+\beta^k)/\da\eq v_k(A,B)/\da\ (\mo\ p)&\t{if}\ (\f B{p^a})=-1.\endcases$$
Also,
$$\beta^k+\l(\f B{p^a}\r)\al^k\eq\cases\al^k+\beta^k\eq v_k(A,B)\ (\mo\ p)&\t{if}\ (\f B{p^a})=1,
\\\beta^k-\al^k\eq -\da u_k(A,B)\ (\mo\ p)&\t{if}\ (\f B{p^a})=-1.\endcases$$
So the desired results follow from the above. \qed

\medskip
\noindent{\it Proof of Theorem 1.1}. Simply apply Theorem 2.1 with $B=m^2$. \qed

\medskip
\noindent{\it Proof of Corollary 1.1}. Let $\omega$ be the primitive cubic root
 $(-1+\sqrt{-3})/2$ of unity. It is easy to see that
$$\l(\f k3\r)=\f{\omega^k-\bar\omega^k}{\sqrt{-3}}=u_k(\omega+\bar\omega,\omega\bar\omega)=u_k(-1,1)$$
for all $k\in\N$. Since $p\eq1\ (\mo\ 3)$, we have $(\f{-3}{p})=(\f {p}3)=1$ and hence
$\da^2\eq-3\ (\mo\ p)$ for some $\da\in\Z$. Observe that
$$\l(\f{-1+\da}{p}\r)^3=\l(\f{\da(\da^2+3)-(3\da^2+1)}{p}\r)=\l(\f{(-3)^2-1}{p}\r)=\l(\f 2{p}\r).$$
By Theorem 1.1,
$$\sum_{k=0}^{p^a-1}\f{u_k(-1,1)\bi{2k}k^h}{(-4)^{hk}}\eq0\ (\mo\ p)\ \ \t{for all}\ h\in\Z^+.$$
If $p^a\eq1\ (\mo\ 12)$, then $(\f{-1}{p^a})=1$ and hence by Theorem 1.1 with $A=m=-1$ we have
$$\sum_{k=0}^{p^a-1}\f{u_k(-1,1)\bi{2k}k^h}{(-1)^k(-4)^{hk}}\eq0\ (\mo\ p)\ \ \t{for all}\ h\in\Z^+.$$

Now assume that $p^a\eq7\ (\mo\ 12)$. Then $(\f{-1+\da}{p^a})=-(\f{-2}{p^a})$ and hence by Theorem 1.1 we have
$$\sum_{k=0}^{p^a-1}\f{v_k(-1,1)\bi{2k}k^h}{(-1)^k(-4)^{hk}}\eq0\ (\mo\ p)\ \ \t{for all}\ h\in\Z^+.$$
Note that
$$v_k(-1,1)=\omega^k+\bar\omega^k=\cases2&\t{if}\ 3\mid k,\\-1&\t{if}\ 3\nmid k.\endcases$$
Thus
$$3\sum^{p^a-1}\Sb k=0\\3\mid k\endSb\f{\bi{2k}k}{(-1)^k(-4)^k}\eq\sum_{k=0}^{p^a-1}\f{\bi{2k}k}{4^k}\eq0\ (\mo\ p).$$
(We apply [ST2, Corollary 1.1] in the last step.) Also,
$$3\sum^{p^a-1}\Sb k=0\\3\mid k\endSb\f{\bi{2k}k^2}{(-1)^k(-4)^{2k}}
\eq\sum_{k=0}^{p^a-1}\f{\bi{2k}k^2}{(-16)^k}\eq\sum_{k=0}^{n}(-1)^k\bi nk^2\ (\mo\ p).$$
where $n=(p^a-1)/2$. Note that
$$\sum_{k=0}^n(-1)^k\bi nk^2=\sum_{k=0}^n(-1)^k\bi nk\bi n{n-k}$$
coincides with the coefficient of $x^n$ in $(1-x)^n(1+x)^n=(1-x^2)^n$ which is zero since $n$ is odd.
Therefore we also have the last two congruences in Corollary 1.1.

The proof of Corollary 1.1 is now complete. \qed

\medskip
\noindent{\it Proof of Corollary 1.2}. As $p\eq\pm1\ (\mo\ 5)$, we have $(\f 5p)=(\f p5)=1$. Thus
$\da^2\eq 5\ (\mo\ p)$ for some $\da\in\Z$. Note that
$$\l(\f 2p\r)\l(\f{3+\da}p\r)=\l(\f{6+2\da}p\r)=\l(\f{(1+\da)^2}p\r)=1.$$
Since $u_k(3,1)=F_{2k}$ and $v_k(3,1)=L_{2k}$, applying Theorem 1.1 with $A=3$ and $m=\pm1$ we immediately obtain the desired results.
\qed

\heading{3. Proofs of Theorems 1.2-1.3 and Corollaries 1.3-1.6}\endheading

\proclaim{Lemma 3.1} Let $p$ be an odd prime and let $x$ be any
algebraic $p$-adic integer. Then
$$\sum_{k=0}^{p-1}\f{\bi{2k}k^2}{16^k}\l(x^k-(-1)^{(p-1)/2}(1-x)^k\r)\eq0\
(\mo\ p^2).$$
\endproclaim
\Proof. This is a result recently obtained by Zhi-Hong Sun [S2] and R. Tauraso [T] independently. \qed

\medskip
\noindent{\it Proof of Theorem 1.2}. (i) Let $\al$ and $\beta$ be the
two roots of the equation $x^2-Ax+B=0$. Then
$$\align&(-1)^{(p-1)/2}\sum_{k=0}^{p-1}\f{\bi{2k}k^2}{(16A)^k}v_k(A,B)
\\=&(-1)^{(p-1)/2}\sum_{k=0}^{p-1}\f{\bi{2k}k^2}{16^k}\(\l(\f{\al}A\r)^k+\l(\f{\beta}A\r)^k\)
\\\eq&\sum_{k=0}^{p-1}\f{\bi{2k}k^2}{16^k}\(\l(1-\f{\al}A\r)^k+\l(1-\f{\beta}A\r)^k\)
\\\eq&\sum_{k=0}^{p-1}\f{\bi{2k}k^2}{16^k}\(\l(\f{\beta}A\r)^k+\l(\f{\al}A\r)^k\)
=\sum_{k=0}^{p-1}\f{\bi{2k}k^2}{(16A)^k}v_k(A,B)\ (\mo\ p^2).
\endalign$$
Hence $$\sum_{k=0}^{p-1}\f{\bi{2k}k^2}{(16A)^k}v_k(A,B)\eq0\ (\mo\
p^2)$$ if $p\eq3\ (\mo\ 4)$.
 Similarly,
$$\align&(-1)^{(p-1)/2}(\al-\beta)\sum_{k=0}^{p-1}\f{\bi{2k}k^2}{(16A)^k}u_k(A,B)
\\=&(-1)^{(p-1)/2}\sum_{k=0}^{p-1}\f{\bi{2k}k^2}{16^k}\(\l(\f{\al}A\r)^k-\l(\f{\beta}A\r)^k\)
\\\eq&\sum_{k=0}^{p-1}\f{\bi{2k}k^2}{16^k}\(\l(1-\f{\al}A\r)^k-\l(1-\f{\beta}A\r)^k\)
\\\eq&\sum_{k=0}^{p-1}\f{\bi{2k}k^2}{16^k}\(\l(\f{\beta}A\r)^k-\l(\f{\al}A\r)^k\)
=(\beta-\al)\sum_{k=0}^{p-1}\f{\bi{2k}k^2}{(16A)^k}u_k(A,B)\ (\mo\
p^2).
\endalign$$
If $p\eq1\ (\mo\ 4)$ and $\Delta=(\al-\beta)^2\not\eq0\ (\mo\ p)$,
then
$$\sum_{k=0}^{p-1}\f{\bi{2k}k^2}{(16A)^k}u_k(A,B)\eq0\ (\mo\ p^2).$$
 So part (i)  holds.

(ii) Below we assume that $(\f{\Delta}p)=1$. Choose $\da\in\Z$ such that
$\da^2\eq\Delta\ (\mo\ p)$. Combining part (i) with Theorem 2.1 in the case $m=A$ and $h=2$, we obtain the
second part of Theorem 1.2.

The proof of Theorem 1.2 is now complete. \qed

\medskip\noindent{\it Proofs of Corollaries 1.3-1.6}.
Recall that
$$\l(\f k3\r)=u_k(-1,1),\ \ F_{2k}=u_k(3,1),\ \ L_{2k}=v_k(3,1),$$
and
$$P_k=u_k(2,-1),\ Q_k=v_k(2,-1),
\ S_k=u_k(4,1)\ T_k=v_k(4,1).$$
In view of this, we immediately obtain the desired results from Theorem 1.2. \qed

\proclaim{Lemma 3.2} Let $A,B\in\Z$. Let $p$ be an odd prime with $(\f Bp)=1$. Suppose that $b^2\eq B\ (\mo\ p)$ where $b\in\Z$.
 Then
$$u_{(p-1)/2}(A,B)\eq\cases0\ (\mo\ p)&\t{if}\ (\f{A^2-4B}p)=1,\\\f1b(\f{A-2b}p)\ (\mo\ p)&\t{if}\ (\f{A^2-4B}p)=-1;
\endcases$$
$$u_{(p+1)/2}(A,B)\eq\cases(\f{A-2b}p)\ (\mo\ p)&\t{if}\ (\f{A^2-4B}p)=1,
\\0\ (\mo\ p)&\t{if}\ (\f{A^2-4B}p)=-1.
\endcases$$
Also,
$$v_{(p-1)/2}(A,B)\eq\cases2(\f{A-2b}p)\ (\mo\ p)&\t{if}\ (\f{A^2-4B}p)=1,\\\ -\f Ab(\f{A-2b}p)
\ (\mo\ p)&\t{if}\ (\f{A^2-4B}p)=-1.\endcases$$
\endproclaim
\Proof. The first two congruences are known results, see, e.g., [S1].
The last one follows from the first two since
$v_n=2u_{n+1}-Au_n$ for $n\in\N$. \qed

\medskip
\noindent{\it Proof of Theorem 1.3}. Let $n=(p-1)/2$, and
$$\al=\f{A+\sqrt{\Delta}}2\ \ \t{and}\ \ \beta=\f{A-\sqrt{\Delta}}2.$$
By Lemma 2.2,
$$\bi {2k}k\eq\bi nk{(-4)^k}\ (\mo\ p)\quad\t{for all}\ k=0,\ldots,p-1.$$
So we have
$$\align&(\al-\beta)\sum_{k=0}^{p-1}\f{u_k(A,B)}{m^k}\bi{2k}k
\\\eq&\sum_{k=0}^{n}\bi nk\(\f{{(-4\al)}^k}{m^k}-\f{{(-4\beta)}^k}{m^k}\)
=\l(1-\f{4\al}m\r)^n-\l(1-\f{4\beta}m\r)^n
\endalign$$
and
$$\align&\sum_{k=0}^{p-1}\f{v_k(A,B)}{m^k}\bi{2k}k
\\\eq&\sum_{k=0}^{n}\bi nk\(\f{{(-4\al)}^k}{m^k}+\f{{(-4\beta)}^k}{m^k}\)
=\l(1-\f{4\al}m\r)^n+\l(1-\f{4\beta}m\r)^n.
\endalign$$
Observe that
$$(m-4\al)+(m-4\beta)=2m-4A\ \t{and}\ (m-4\al)(m-4\beta)=m^2-4mA+16B.$$
Thus
$$\align\l(\f mp\r)\sum_{k=0}^{p-1}\f{u_k(A,B)}{m^k}\bi{2k}k\eq&-4\times\f{(m-4\al)^n-(m-4\beta)^n}{4\beta-4\al}
\\\eq&-4u_n(2m-4A,m^2-4Am+16B)\ (\mo\ p)\endalign$$
and
$$\l(\f mp\r)\sum_{k=0}^{p-1}\f{v_k(A,B)}{m^k}\bi{2k}k\eq v_n(2m-4A,m^2-4Am+16B)\ (\mo\ p).$$
Note that
$$(2m-4A)^2-4(m^2-4Am+16B)=16(A^2-4B)=16\Delta.$$
Via Lemma 3.2 we are able to determine $u_n(2m-4A,m^2-4Am+16B)$ and $v_n(2m-4A,m^2-4Am+16B)$
modulo $p$ and hence the desired congruences follow. \qed

\heading{4. Proofs of Theorems 1.4 and 1.5}\endheading

\medskip{\it Proof of Theorem 1.4}. As $(\f{\Delta}p)=1$, there is an integer $\da$ such that $\da^2\eq\Delta\ (\mo\ p^2)$.
Set $\al=(A+\da)/2$ and $\beta=(A-\da)/2$. Then
$$\align&\sum_{k=0}^{p-1}\f{A^kv_k(A,B)}{(4B)^k}\bi{2k}k
\\\eq&\sum_{k=0}^{p-1}\bi{2k}k\(\f{(A\al)^k}{(4B)^k}+\f{(A\beta)^k}{(4B)^k}\)
=\sum_{k=0}^{p-1}\f{\bi{2k}k}{(4\beta/A)^k}+\sum_{k=0}^{p-1}\f{\bi{2k}k}{(4\al/A)^k}\ (\mo\ p^2).
\endalign$$
Note that
$$\f{4\al}A\l(\f{4\al}A-4\r)\eq-\f{4^2}{A^2}B\eq\f{4\beta}A\l(\f{4\beta}A-4\r)\ (\mo\ p^2).$$
Hence by the main result of [S09b] we have
$$\align &\sum_{k=0}^{p-1}\f{\bi{2k}k}{(4\al/A)^k}+\sum_{k=0}^{p-1}\f{\bi{2k}k}{(4\beta/A)^k}
\\\eq&\l(\f{-B}p\r)+u_{p-(\f{-B}p)}\l(\f{4\al}A-2,1\r)
+\l(\f{-B}p\r)+u_{p-(\f{-B}p)}\l(\f{4\beta}A-2,1\r)
\ (\mo\ p^2).
\endalign$$
Since $$\f{4\al}A-2+\f{4\beta}A-2=0$$
and
$u_n(-x,1)=(-1)^{n-1}u_n(x,1)$ for any $n\in\N$, the desired result follows from the above. \qed

\proclaim{Lemma 4.1} Let $p\not=2,5$ be a prime. Then
$$F_{(p-(\f p5))/2}\eq\cases0\ (\mo\ p)&\t{if}\ p\eq1\ (\mo\ 4),
\\2(-1)^{\lfloor(p+5)/10\rfloor}(\f 5p)5^{(p-3)/4}\ (\mo\ p)&\t{if}\ p\eq 3\ (\mo\ 4),
\endcases$$
and
$$F_{(p+(\f p5))/2}\eq\cases(-1)^{\lfloor(p+5)/10\rfloor}(\f 5p)5^{(p-1)/4}\ (\mo\ p)&\t{if}\ p\eq1\ (\mo\ 4),
\\(-1)^{\lfloor(p+5)/10\rfloor}(\f 5p)5^{(p-3)/4}\ (\mo\ p)&\t{if}\ p\eq 3\ (\mo\ 4).
\endcases$$
Also,
$$L_{(p-(\f p5))/2}\eq\cases2(-1)^{\lfloor(p+5)/10\rfloor}(\f 5p)5^{(p-1)/4}\ (\mo\ p)&\t{if}\ p\eq1\ (\mo\ 4),
\\0\ (\mo\ p)&\t{if}\ p\eq 3\ (\mo\ 4),
\endcases$$
and
$$L_{(p+(\f p5))/2}\eq\cases(-1)^{\lfloor(p+5)/10\rfloor}5^{(p-1)/4}\ (\mo\ p)&\t{if}\ p\eq1\ (\mo\ 4),
\\(-1)^{\lfloor(p+5)/10\rfloor}(\f 5p)5^{(p+1)/4}\ (\mo\ p)&\t{if}\ p\eq 3\ (\mo\ 4).
\endcases$$
\endproclaim
\Proof. This follows from Z. H. Sun and Z. W. Sun [SS, Corollaries 1 and 2]. \qed

\medskip
\noindent{\it Proof of Theorem 1.5}.  As in the proof of Theorem 1.3, we have
$$\align\sum_{k=0}^{p-1}\f{F_k}{12^k}\bi{2k}k\eq&-4\l(\f{12}p\r)u_{(p-1)/2}(2\times 12-4,12^2-4\times1\times12+16(-1))
\\\eq&-4\l(\f 3p\r)u_{(p-1)/2}(20,80)\ (\mo\ p).
\endalign$$
Set $n=(p-1)/2$. As the equations $x^2-20x+80=0$ has two roots $10\pm2\sqrt5$, we have
$$\align u_n(20,80)=&\f{(10+2\sqrt5)^n-(10-2\sqrt5)^n}{4\sqrt5}
\\=&(4\sqrt5)^{n-1}\(\(\f{1+\sqrt5}2\)^n-(-1)^n\(\f{1-\sqrt5}2\)^n\)
\\=&\cases2^{p-3}5^{(p-1)/4}F_n&\t{if}\ 2\mid n,\ \t{i.e.,}\ p\eq 1\ (\mo\ 4),
\\ 2^{p-3}5^{(p-3)/4}L_n&\t{if}\ 2\nmid n, \ \t{i.e.,}\ p\eq 3\ (\mo\ 4).\endcases
\endalign$$
Therefore
$$-\l(\f 3p\r)\sum_{k=0}^{p-1}\f{F_k}{12^k}\bi{2k}k\eq\cases 5^{(p-1)/4}F_{(p-1)/2}&\t{if}\ p\eq1\ (\mo\ 4),
\\5^{(p-3)/4}L_{(p-1)/2}&\t{if}\ p\eq 3\ (\mo\ 4).\endcases$$

{\it Case}\ 1. $(\f p5)=1$. By Lemma 4.1,
$$F_{(p-1)/2}=F_{(p-(\f p5))/2}\eq0\ (\mo\ p)\quad\t{if}\ p\eq1\ (\mo\ 4),$$
and
$$L_{(p-1)/2}=L_{(p-(\f p5))/2}\eq0\ (\mo\ p)\quad\t{if}\ p\eq3\ (\mo\ 4).$$
It follows that
$$\sum_{k=0}^{p-1}\f{F_k}{12^k}\bi{2k}k\eq0\ (\mo\ p).$$

{\it Case}\ 2. $(\f p5)=-1$. If $p\eq 1\ (\mo\ 4)$, then by Lemma 4.1 we have
$$\align&5^{(p-1)/4}F_{(p-1)/2}=5^{(p-1)/4}F_{(p+(\f p5))/2}
\\\eq&5^{(p-1)/4}(-1)^{\lfloor(p+5)/10\rfloor}\l(\f 5p\r)5^{(p-1)/4}\eq(-1)^{\lfloor(p+5)/10\rfloor} \ (\mo\ p).
\endalign$$
If $p\eq 3\ (\mo\ 4)$, then by Lemma 4.1 we have
$$\align&5^{(p-3)/4}L_{(p-1)/2}=5^{(p-3)/4}L_{(p+(\f p5))/2}
\\\eq&5^{(p-3)/4}(-1)^{\lfloor(p+5)/10\rfloor}\l(\f 5p\r)5^{(p+1)/4}\eq(-1)^{\lfloor(p+5)/10\rfloor} \ (\mo\ p).
\endalign$$
Therefore
$$\sum_{k=0}^{p-1}\f{F_k}{12^k}\bi{2k}k\eq-\l(\f p3\r)(-1)^{\lfloor(p+5)/10\rfloor}\ (\mo\ p)$$
and hence the first congruence in Theorem 1.5 follows.

 The second congruence in Theorem 1.5 can be proved in a similar way. We omit the details. \qed

\heading{5. Some conjectures}\endheading

 Our following conjectures involve representations of primes by binary quadratic forms.
 The reader may consult [C] and [BEW, Chapter 9] for basic knowledge and background.

\proclaim{Conjecture 5.1} Let $p>3$ be a prime.
If $p\eq7\ (\mo\ 12)$ and $p=x^2+3y^2$ with $y\eq1\ (\mo\ 4)$, then
$$\sum_{k=0}^{p-1}\l(\f k3\r)\f{\bi{2k}k^2}{(-16)^k}\eq (-1)^{(p-3)/4}\l(4y-\f p{3y}\r)\ (\mo\ p^2)$$
and
$$\sum_{k=0}^{p-1}\l(\f k3\r)\f{k\bi{2k}k^2}{(-16)^k}\eq (-1)^{(p+1)/4}y\ (\mo\ p^2).$$
If $p\eq11\ (\mo\ 12)$, then
$$\sum_{k=0}^{p-1}\l(\f k3\r)\f{\bi{2k}k^2}{(-16)^k}\eq0\ (\mo\ p).$$
If $p\eq1\ (\mo\ 12)$, then
$$\sum_{k=0}^{p-1}\bi{p-1}k\l(\f k3\r)\f{\bi{2k}k^2}{16^k}\eq0\ (\mo\ p^2).$$
\endproclaim

\proclaim{Conjecture 5.2} {\rm (i)} Let $p$ be a prime
with $p\eq1,3\ (\mo\ 8)$. Write $p=x^2+2y^2$ with $x,y\in\Z$ and $x\eq1,3\ (\mo\ 8)$. Then
$$\sum_{k=0}^{p-1}\f{P_k}{(-8)^k}\bi{2k}k^2\eq \cases0\ (\mo\ p^2)&\t{if}\ p\eq1\ (\mo\ 8),
\\(-1)^{(p-3)/8}(p/(2x)-2x)\ (\mo\ p^2)&\t{if}\ p\eq3\ (\mo\ 8).
\endcases$$
Also,
$$\sum_{k=0}^{p-1}\f{kP_k}{(-8)^k}\bi{2k}k^2\eq\f{(-1)^{(x+1)/2}}2\l(x+\f p{2x}\r)\ (\mo\ p^2).$$

{\rm (ii)} If $p\eq5\ (\mo\ 8)$ is a prime, then
$$\sum_{k=0}^{p-1}\f{P_k}{(-8)^k}\bi{2k}k^2\eq0\ (\mo\ p).$$
If $p\eq7\ (\mo\ 8)$ is a prime, then
$$\sum_{k=0}^{p-1}\bi{p-1}k\f{P_k}{8^k}\bi{2k}k^2\eq0\ (\mo\ p^2).$$
\endproclaim

\proclaim{Conjecture 5.3} Let $p$ be an odd prime.

{\rm (i)} If $p\eq3\ (\mo\ 8)$ and $p=x^2+2y^2$ with $y\eq1,3\ (\mo\ p)$, then
$$\sum_{k=0}^{p-1}\f{P_k}{32^k}\bi{2k}k^2\eq (-1)^{(y-1)/2}\l(2y-\f p{4y}\r)\ (\mo\ p^2).$$
If $p\eq 7\ (\mo\ 8)$, then
$$\sum_{k=0}^{p-1}\f{P_k}{32^k}\bi{2k}k^2\eq0\ (\mo\ p).$$

{\rm (ii)} Suppose that $p\eq1,3\ (\mo\ 8)$, $p=x^2+2y^2$ with $x\eq1,3\ (\mo\ 8)$ and also $y\eq1,3\ (\mo\ 8)$ when $p\eq3\ (\mo\ 8)$.
Then
$$\sum_{k=0}^{p-1}\f{kP_k}{32^k}\bi{2k}k^2\eq\cases(-1)^{(p-1)/8}(p/(4x)-x/2)\ (\mo\ p^2)&\t{if}\ p\eq1\ (\mo\ 8),
\\(-1)^{(y+1)/2}y\ (\mo\ p^2)&\t{if}\ p\eq3\ (\mo\ 8).
\endcases$$
\endproclaim

\proclaim{Conjecture 5.4} Let $p$ be an odd prime.

{\rm (i)} When $p\eq1,3\ (\mo\ 8)$ and $p=x^2+2y^2$ with $x,y\in\Z$ and $x\eq1,3\ (\mo\ 8)$, we have
$$\sum_{k=0}^{p-1}\f{Q_k}{(-8)^k}\bi{2k}k^2\eq (-1)^{(x-1)/2}\l(4x-\f px\r)\ (\mo\ p^2)$$
and
$$\sum_{k=0}^{p-1}\f{kQ_k}{(-8)^k}\bi{2k}k^2\eq\cases0\ (\mo\ p^2)&\t{if}\ p\eq1\ (\mo\ 8),
\\(-1)^{(p-3)/8}2(x+p/x)\ (\mo\ p^2)&\t{if}\ p\eq3\ (\mo\ 8).
\endcases$$

{\rm (ii)} When $p\eq5,7\ (\mo\ 8)$, we have
$$\sum_{k=0}^{p-1}\f{Q_k}{(-8)^k}\bi{2k}k^2\eq0\ (\mo\ p).$$
\endproclaim

\proclaim{Conjecture 5.5} Let $p$ be an odd prime.

{\rm (i)} When $p\eq1\ (\mo\ 8)$ and $p=x^2+2y^2$ with $x,y\in\Z$ and $x\eq1,3\ (\mo\ 8)$, we have
$$\sum_{k=0}^{p-1}\f{Q_k}{32^k}\bi{2k}k^2\eq (-1)^{(p-1)/8}\l(4x-\f px\r)\ (\mo\ p^2).$$
If $p\eq 5\ (\mo\ 8)$, then
$$\sum_{k=0}^{p-1}\f{Q_k}{32^k}\bi{2k}k^2\eq 0\ (\mo\ p).$$

{\rm (ii)} If $p\eq1,3\ (\mo\ 8)$ and $p=x^2+2y^2$ with $x\eq1,3\ (\mo\ 8)$ and also $y\eq1,3\ (\mo\ 8)$ when $p\eq3\ (\mo\ 8)$,
then
$$\sum_{k=0}^{p-1}\f{k\bi{2k}k^2}{32^k}Q_k\eq\cases(-1)^{(p-1)/8}(p/x-2x)\ (\mo\ p^2)&\t{if}\ p\eq1\ (\mo\ 8),
\\(-1)^{(y+1)/2}2y\ (\mo\ p^2)&\t{if}\ p\eq 3\ (\mo\ 8).\endcases$$
\endproclaim

\proclaim{Conjecture 5.6} Let $p>3$ be a prime.
If $p\eq7\ (\mo\ 12)$ and $p=x^2+3y^2$ with $y\eq1\ (\mo\ 4)$, then
$$\sum_{k=0}^{p-1}\f{S_k}{4^k}\bi{2k}k^2\eq (-1)^{(p+1)/4}\l(4y-\f p{3y}\r)\ (\mo\ p^2)$$
and
$$\sum_{k=0}^{p-1}\f{kS_k}{4^k}\bi{2k}k^2\eq(-1)^{(p-3)/4}\l(6y-\f{7p}{3y}\r)\ (\mo\ p^2).$$
We also have
$$\sum_{k=0}^{p-1}\f{S_k}{4^k}\bi{2k}k^2\eq\cases0\ (\mo\ p^2)&\t{if}\ p\eq1\ (\mo\ 12),
\\0\ (\mo\ p)&\t{if}\ p\eq2\ (\mo\ 3).\endcases$$
\endproclaim

\proclaim{Conjecture 5.7} Let $p>3$ be a prime.
If $p\eq7\ (\mo\ 12)$ and $p=x^2+3y^2$ with $y\eq1\ (\mo\ 4)$, then
$$\sum_{k=0}^{p-1}\f{S_k}{64^k}\bi{2k}k^2\eq 2y-\f p{6y}\ (\mo\ p^2)$$
and
$$\sum_{k=0}^{p-1}\f{kS_k}{64^k}\bi{2k}k^2\eq y\ (\mo\ p^2).$$
If $p\eq 11\ (\mo\ 12)$, then
$$\sum_{k=0}^{p-1}\f{S_k}{64^k}\bi{2k}k^2\eq0\ (\mo\ p).$$
\endproclaim

\proclaim{Conjecture 5.8} Let $p>3$ be a prime.

{\rm (i)} If $p\eq1\ (\mo\ 12)$ and $p=x^2+3y^2$ with $x\eq1\ (\mo\ 3)$, then
$$\sum_{k=0}^{p-1}\f{T_k}{4^k}\bi{2k}k^2\eq(-1)^{(p-1)/4+(x-1)/2}\l(4x-\f p{x}\r)\ (\mo\ p^2)$$
and
$$\sum_{k=0}^{p-1}\f{T_k}{64^k}\bi{2k}k^2\eq(-1)^{(x-1)/2}\l(4x-\f p{x}\r)\ (\mo\ p^2);$$
also
$$\sum_{k=0}^{p-1}\f{kT_k}{4^k}\bi{2k}k^2\eq (-1)^{(p-1)/4+(x+1)/2}\l(4x-\f{2p}x\r)\ (\mo\ p^2)$$
and
$$\sum_{k=0}^{p-1}\f{kT_k}{64^k}\bi{2k}k^2\eq (-1)^{(x-1)/2}\l(2x-\f{p}x\r)\ (\mo\ p^2).$$

{\rm (ii)} If $p\eq7\ (\mo\ 12)$ and $p=x^2+3y^2$ with $y\eq1\ (\mo\ 4)$, then
$$\sum_{k=0}^{p-1}\f{T_k}{4^k}\bi{2k}k^2\eq (-1)^{(p-3)/4}\l(12y-\f p{y}\r)\ (\mo\ p^2),$$
$$\sum_{k=0}^{p-1}\f{kT_k}{4^k}\bi{2k}k^2\eq(-1)^{(p+1)/4}\l(20y-\f{8p}y\r)\ (\mo\ p^2)$$
and
$$\sum_{k=0}^{p-1}\f{kT_k}{64^k}\bi{2k}k^2\eq 4y\ (\mo\ p^2).$$

{\rm (iii)} If $p\eq 5\ (\mo\ 12)$, then
$$\sum_{k=0}^{p-1}\f{T_k}{4^k}\bi{2k}k^2\eq\sum_{k=0}^{p-1}\f{T_k}{64^k}\bi{2k}k^2\eq0\ (\mo\ p).$$
If $p\eq 11\ (\mo\ 12)$, then
$$\sum_{k=0}^{p-1}\bi{p-1}k\f{T_k}{(-4)^k}\bi{2k}k^2\eq0\ (\mo\ p^2).$$
\endproclaim


\medskip

 \widestnumber\key{BEW}

 \Refs

\ref\key BEW\by B. C. Berndt, R. J. Evans and K. S. Williams
\book Gauss and Jacobi Sums\publ John Wiley \& Sons, 1998\endref

\ref\key C\by D. A. Cox\book Primes of the Form $x^2+ny^2$\publ John Wiley \& Sons, 1989\endref

\ref\key M1\by E. Mortenson\paper A supercongruence conjecture of Rodriguez-Villegas
for a certain truncated hypergeometric function
\jour J. Number Theory\vol 99\yr 2003\pages 139--147\endref


\ref\key M2\by E. Mortenson\paper Supercongruences between truncated ${}_2\! F_1$
by geometric functions and their Gaussian analogs
\jour Trans. Amer. Math. Soc.\vol 355\yr 2003\pages 987--1007\endref

\ref\key M3\by E. Mortenson\paper Supercongruences for truncated  ${}_{n+1}\! F_n$
hypergeometric series with applications to certain weight three newforms
\jour Proc. Amer. Math. Soc.\vol 133\yr 2005\pages 321--330\endref

\ref\key O\by K. Ono\book Web of Modularity: Arithmetic of the Coefficients of Modular Forms and $q$-series
\publ Amer. Math. Soc., Providence, R.I., 2003\endref

\ref\key PS\by H. Pan and Z. W. Sun\paper A combinatorial identity
with application to Catalan numbers \jour Discrete Math.\vol
306\yr 2006\pages 1921--1940\endref


\ref\key RV\by F. Rodriguez-Villegas\paper Hypergeometric families of Calabi-Yau manifolds
\jour {\rm in}: Calabi-Yau Varieties and Mirror Symmetry (Toronto, ON, 2001), pp. 223-231,
Fields Inst. Commun., {\bf 38}, Amer. Math. Soc., Providence, RI, 2003\endref

\ref\key St1\by R. P. Stanley\book Enumerative Combinatorics \publ
Vol. 1, Cambridge Univ. Press, Cambridge, 1999\endref

\ref\key St2\by R. P. Stanley\book Enumerative Combinatorics \publ
Vol. 2, Cambridge Univ. Press, Cambridge, 1999\endref

\ref\key S1\by Z. H. Sun\paper Values of Lucas sequences modulo primes
\jour Rocky Mount. J. Math. \vol 33\yr 2003\pages 1123--1145\endref

\ref\key S2\by Z. H. Sun\paper Congruences concerning Legendre polynomials
\jour preprint, 2009\endref

\ref\key SS\by Z. H. Sun and Z. W. Sun\paper Fibonacci numbers and Fermat's last theorem
\jour Acta Arith.\vol 60\yr 1992\pages 371--388\endref

\ref\key S02\by Z. W. Sun\paper On the sum $\sum_{k\eq r\,(\mo\ m)}\bi nk$
and related congruences\jour Israel J. Math.
\vol 128\yr 2002\pages 135--156\endref


\ref\key S09a\by Z. W. Sun\paper Various congruences involving binomial coefficients and higher-order Catalan numbers
\jour  arXiv:0909.3808. {\tt http://arxiv.org/abs/0909.3808}\endref

\ref\key S09b\by Z. W. Sun\paper Binomial coefficients, Catalan numbers and Lucas quotients
\jour preprint, arXiv:0909.5648. {\tt http://arxiv.org/abs/0909.5648}\endref

\ref\key S09c\by Z. W. Sun\paper $p$-adic valuations of some sums of multinomial coefficients
\jour preprint, arXiv:0910.3892. {\tt http://arxiv.org/abs/0910.3892}\endref

\ref\key S09d\by Z. W. Sun\paper On sums of binomial coefficients modulo $p^2$
\jour preprint, arXiv:0910.5667. {\tt http://arxiv.org/abs/0910.5667}\endref

\ref\key S09e\by Z. W. Sun\paper Binomial coefficients, Catalan numbers and Lucas quotients (II)
\jour preprint, arXiv:0911.3060. {\tt http://arxiv.org/abs/0911.3060}\endref

\ref\key S09f\by Z. W. Sun\paper On congruences related to central binomial coefficients
\jour preprint, arXiv:0911.2415. {\tt http://arxiv.org/abs/0911.2415}\endref

\ref\key ST1\by Z. W. Sun and R. Tauraso\paper On some new congruences for binomial coefficients
\jour Acta Arith.\pages to appear. {\tt http://arxiv.org/abs/0709.1665}\endref

\ref\key ST2\by Z. W. Sun and R. Tauraso\paper New congruences for central binomial coefficients
\jour Adv. in Appl. Math., to appear. {\tt http://arxiv.org/abs/0805.0563}\endref

\ref\key T\by R. Tauraso\paper An elementary proof of a Rodriguez-Villegas supercongruence
\jour preprint, arXiv:0911.4261. {\tt http://arxiv.org/abs/0911.4261}\endref

\endRefs
\enddocument